\documentclass{article}

\usepackage[utf8]{inputenc}
\usepackage{CJKutf8}

\usepackage[inline]{asymptote}

\usepackage[version=4]{mhchem}

\usepackage{array}

\usepackage{graphicx}
\graphicspath{ {./images/} }

\usepackage[usenames,dvipsnames]{color}
\definecolor{dkgreen}{rgb}{0,0.6,0}
\definecolor{gray}{rgb}{0.5,0.5,0.5}
\definecolor{mauve}{rgb}{0.58,0,0.82}
\definecolor{lavender}{rgb}{0.8,0.56,0.98}
\definecolor{emerald}{RGB}{21, 89, 26}
\definecolor{irisElfLipgloss}{RGB}{199, 56, 65}

\usepackage{tikz}
\usetikzlibrary{arrows.meta,shapes.misc,decorations.pathmorphing,calc,bending,decorations.markings,shapes.geometric}

\usepackage[makeroom]{cancel}

\usepackage[normalem]{ulem}

\usepackage{amsmath}
\usepackage{mathtools}
\usepackage{amssymb}
\usepackage{wasysym}
\usepackage{amsfonts}
\usepackage{mathrsfs}
\DeclareMathAlphabet{\mathpzc}{OT1}{pzc}{m}{it}

\usepackage{bigints}

\usepackage{relsize}

\usepackage[english]{babel} 
\usepackage[autostyle, english = american]{csquotes}
\MakeOuterQuote{"} 

\usepackage{hyperref}
\hypersetup{
    colorlinks,
    linkcolor={red!50!black},
    citecolor={blue!50!black},
    urlcolor={blue!80!black}
}

\usepackage{geometry}
\usepackage{mdframed}
\usepackage{calc}
\definecolor{mytitlecolor}{RGB}{0,102,204}
\newenvironment{theorem}[1][]{%
    \begin{mdframed}[
        backgroundcolor=mytitlecolor!6,
        linewidth=0pt,
        innerleftmargin=10pt,
        innerrightmargin=10pt,
        innertopmargin=10pt,
        innerbottommargin=10pt,
        leftmargin=\dimexpr(\linewidth-0.92\linewidth)/2\relax, 
        rightmargin=\dimexpr(\linewidth-0.92\linewidth)/2\relax,
        skipabove=\topsep,
        skipbelow=\topsep
    ]
    \ifx\relax#1\relax\else\textcolor{mytitlecolor}{\large\textbf{#1}:}\par\fi
    }{
    \end{mdframed}
}

\usepackage{tcolorbox}
\definecolor{exercolor}{RGB}{0,200,115}

\usepackage{tcolorbox}
\definecolor{examplecolor}{RGB}{161,161,161}

\usepackage{mdframed}
\usepackage{calc}
\definecolor{myDefinitionColor}{RGB}{255, 104, 0}
\newenvironment{definition}[1][]{%
    \begin{mdframed}[
        backgroundcolor=myDefinitionColor!6,
        linewidth=0pt,
        innerleftmargin=10pt,
        innerrightmargin=10pt,
        innertopmargin=10pt,
        innerbottommargin=10pt,
        leftmargin=\dimexpr(\linewidth-0.92\linewidth)/2\relax, 
        rightmargin=\dimexpr(\linewidth-0.92\linewidth)/2\relax,
        skipabove=\topsep,
        skipbelow=\topsep
    ]
    \ifx\relax#1\relax\else\textcolor{myDefinitionColor}{\large\textbf{#1}:}\par\fi
    }{
    \end{mdframed}
}

\usepackage{listings}
\definecolor{dkgreen}{rgb}{0,0.6,0}
\definecolor{gray}{rgb}{0.5,0.5,0.5}
\definecolor{mauve}{rgb}{0.58,0,0.82}
\usepackage{svg}

\makeatletter
\def\smallunderbrace#1{
  \mathop{
    \vtop{
      \m@th
      \ialign{##\crcr
        $\hfil\displaystyle{#1}\hfil$\crcr
        \noalign{\kern3\p@\nointerlineskip}
        \tiny\upbracefill\crcr
        \noalign{\kern3\p@}
      }
    }
  }\limits
}
\makeatother

\usepackage[d]{esvect}

\renewcommand\vec[1]{\vv{\textbf{#1}}}

\usepackage{biblatex}
\DeclareCiteCommand{\supercite}[\mkbibsuperscript]
  {\iffieldundef{prenote}
     {}
     {\BibliographyWarning{Ignoring prenote argument}}%
   \iffieldundef{postnote}
     {}
     {\BibliographyWarning{Ignoring postnote argument}}}
  {\usebibmacro{citeindex}%
   \bibopenbracket\usebibmacro{cite}\bibclosebracket}
  {\supercitedelim}
  {}
\let\cite=\supercite

\usepackage{etoolbox}
\makeatletter
\patchcmd{\l@section}
    {\hfil}
    {\leaders\hbox{\normalfont$\hbox{.}\mkern \@dotsep mu$}\hfill}
    {}{}
\makeatother
\usepackage{multicol}
\begin{document}

\begin{center}\vspace{-1cm}
    \textbf{\LARGE Closed Form of a Generalized Sinkhorn Limit}\\
    Max Chicky Fang
    ~
\end{center}

\begin{center}
\subsection*{Abstract}

    The Kruithof iterative scaling process, which adjusts matrices to meet target row and column sums, is a longstanding problem that lacks a general closed form for its limit. While Nathanson derived the closed form for the Sinkhorn limit of $2\times 2$ matrices when target row and column sums are 1, and recent work by Rowland and Wu has advanced understanding of Sinkhorn limits for $3\times 3$, and general $n\times m$ matrices through polynomials, a "generalized Sinkhorn limit" (i.e. the original "Kruithof limit", with arbitrary target sums) remains elusive. Here, we derive the closed form for the generalized Sinkhorn limit of $2\times 2$ matrices, and discuss how this approach can be extended to larger matrices. More significantly, we prove that for any positive $n \times m$ matrix and positive target row and column sums, each entry in the generalized Sinkhorn limit is algebraic over the input data with degree at most $\binom{n+m-2}{n-1}$.
    
\end{center}

\section{Introduction}
\subsection{Background}

    Iterative scaling procedures for matrices have a long history in mathematics and statistics. A notable example is the Sinkhorn-Knopp algorithm, introduced by Sinkhorn in 1964,\cite{Sinkhorn-1964}\cite{Sinkhorn-Knopp-algorithm-1967} which iteratively scales the rows and columns of a positive square matrix to produce a doubly stochastic matrix, not dissimilar to the iterative proportional fitting procedure in statistics. Sinkhorn showed that for any positive square matrix $A$, this iterative process converges to a unique doubly stochastic matrix $S$, known as the Sinkhorn limit of $A$, which can be expressed as $S=\operatorname{Sink}(A)=D_1AD_2$ where $D_{1,2}$ are positive diagonal matrices, unique up to a scalar multiple. In a subsequent 1966 paper,\cite{Sinkhorn-1966} Sinkhorn extended this result, showing that a single diagonal scaling of the form $DAD$, where $A$ is a positive square matrix and $D$ is a positive diagonal matrix, unique up to a scalar multiple, suffices to produce a singly row or column stochastic matrix.\\
    
    The structure of Sinkhorn limits are now well understood for small matrices, including a closed form of the $2\times2$ case due to Nathanson,\cite{Nathanson-2019} 
    $$
    A=\begin{bmatrix}
        a&b\\
        c&d
    \end{bmatrix}
    \iff
    \operatorname{Sink}(A)=
    \frac{1}{\sqrt{a d}+\sqrt{b c}}
    \begin{bmatrix}
        \sqrt{a d} & \sqrt{b c} \\
        \sqrt{b c} & \sqrt{a d} \\
    \end{bmatrix}
    $$
    and recent work by Rowland and Wu has extended these insights to larger matrices.\cite{rowland-wu-2024} In particular, they derive the algebraic structure of the Sinkhorn limit for general $3\times 3$ matrices using Gröbner bases, proving that each entry is an algebraic number of degree at most 6. Furthermore, they develop a conjectural framework for the Sinkhorn limits of general square $n\times n$, and later, rectangular $n\times m$ matrices, characterizing the entries as roots of polynomials with degree at most $\binom{m+n-2}{n-1}$, whose coefficients are structured as linear combinations of products of matrix minors.\\

    This procedure can be viewed as a special case of a more general scaling method, often referred to as the Kruithof scaling process,\cite{Kruithof-Original-1937} originally introduced in 1937 to model and adjust traffic flows in telephone networks. The goal here was to adjust a general rectangular $n\times m$ matrix to meet prescribed row and column sums, not necessarily equal to $1$. Unlike the Sinkhorn-Knopp case, this "Kruithof limit", or hereinafter \textit{generalized Sinkhorn limit} lacks a closed form for the limiting matrix, even in small cases. While Sinkhorn showed that the limit exists and is unique,\cite{Sinkhorn-proof-of-KruithoffExistence-1967} and Franklin and Lorenz implicitly demonstrate its continuity among other impressive general results,\cite{franklin-verygeneral-1989scaling} much, if not all, of remaining theoretical focus has been on the simpler setting of scaling to the doubly stochastic matrix.\\

    Here, we derive the closed form for the generalized Sinkhorn limits of $2\times2$ matrices with arbitrary positive target row and column sums. Building on earlier results by Sinkhorn, Franklin and Lorenz, Nathanson, Rowland and Wu, \&c., who primarily address the doubly stochastic case, our approach extends these ideas to the arbitrary scaling setting as originally envisioned by Kruithof. Although this is tractable for $2\times 2$ matrices, we find that the algebra quickly becomes unwieldy in larger cases, making explicit formulas for larger matrices likely impractical. 

\subsection{Observations}

    Despite the apparent simplicity of scaling a matrix to meet prescribed row and column sums, the exact values of the entries of Sinkhorn limits or generalized Sinkhorn limits remain remarkably opaque, even in small cases. This difficulty reflects the general unwieldiness of matrix operations and, for instance, their interaction with spectral properties. In operator theory and matrix analysis, it is well known that even basic operations like matrix addition or multiplication can obscure or completely alter eigenstructure. For example, for matrices $A$ and $B$, while the nonzero eigenvalues of matrix products $AB$ and $BA$ coincide up to multiplicity, there is no general formula that relates the eigenvalues of $A+B$ to those of $A$ and $B$ individually. Generally, in finite dimensions, matrix behavior resists clean decomposition. As a result, even when scaling a matrix iteratively, one may encounter nontrivial dynamics that make the algebraic characterization of the limit matrix highly nonobvious.\cite{norberDiscussionAA} The work of Rowland and Wu reinforces this notion through their analysis of Sinkhorn limits, highlighting the nontrivial polynomial structure underlying even small matrices.

\section{Calculation}
\subsection{Setup}

    Let $A\in \mathbb R^{n\times m}_{>0}$ be a positive $n\times m$ rectangular matrix. Denote $S=\operatorname{Sink}\left(A\left|\, \vec{R},\vec{C} \right.\right)$ to be the generalized Sinkhorn limit of $A$ with target row and column sums $\vec{R}\in \mathbb R^{n}_{>0}$ and $\vec{C}\in \mathbb R^{m}_{>0}$ defined respectively as
    $$
    \vec{R} = S \vec{1}_m = \begin{bmatrix} R_1\\  R_2 \\ \vdots \\ R_n \end{bmatrix},\quad 
    \vec{C} = S^\top \vec{1}_n = \begin{bmatrix} C_1 \\ C_2 \\ \vdots \\ C_m \end{bmatrix}
    $$
    for $\vec{1}_i$ being a column vector of $i$ ones. Sinkhorn showed that like the Sinkhorn limit, the generalized Sinkhorn limit can be equivalently written as the product of three matrices $\operatorname{Sink}\left(A\left|\, \vec{R},\vec{C} \right.\right)=D_1 AD_2$ for positive diagonal matrices $D_{1,2}$, unique up to a scalar multiple.\cite{Sinkhorn-proof-of-KruithoffExistence-1967} \\
    
    This means that over successive iterations, the scaling factors applied to each row and column converge, with row and column scalings approaching fixed values $r_i$ and $c_j$ respectively. We observe that $D_1$ can be interpreted as scaling factors $r_i$ for the corresponding $i$th row of $A$, and $D_2$ being scaling factors $c_j$ for the corresponding $j$th column of $A$, 
    $$D_1 AD_2=
    \begin{bmatrix}
        r_1 & 0 & \cdots & 0 \\
        0 & r_2 & \cdots & 0 \\
        \vdots & \vdots & \ddots & \vdots \\
        0 & 0 & \cdots & r_n 
    \end{bmatrix}
    \begin{bmatrix}
        a_{11} & a_{12} & \cdots & a_{1m} \\
        a_{21} & a_{22} & \cdots & a_{2m} \\
        \vdots & \vdots & \ddots & \vdots \\
        a_{n1} & a_{n2} & \cdots & a_{nm}
    \end{bmatrix}
    \begin{bmatrix}
        c_1 & 0 & \cdots & 0 \\
        0 & c_2 & \cdots & 0 \\
        \vdots & \vdots & \ddots & \vdots \\
        0 & 0 & \cdots & c_m 
    \end{bmatrix}
    =
    \begin{bmatrix}
        r_1 a_{11} c_1 & r_1 a_{12} c_2 & \cdots & r_1 a_{1m} c_m \\
        r_2 a_{21} c_1 & r_2 a_{22} c_2 & \cdots & r_2 a_{2m} c_m \\
        \vdots & \vdots & \ddots & \vdots \\
        r_n a_{n1} c_1 & r_n a_{n2} c_2 & \cdots & r_n a_{nm} c_m \\
    \end{bmatrix}
    $$
    Effectively, the generalized Sinkhorn limit converges to a matrix where each entry $a_{ij}$ is transformed into $r_ia_{ij}c_j$. Formally, we can say that
    \begin{definition}[Definition 2.1.1 \normalfont{(Generalized Sinkhorn Limit)}]

        Let $A \in \mathbb{R}^{n \times m}_{>0}$ be a positive $n\times m$ rectangular matrix, and let target row and column sums be $\vec{R} \in \mathbb{R}^n_{>0}$ and $\vec{C} \in \mathbb{R}^m_{>0}$, respectively.\\

        The \textit{generalized Sinkhorn limit} $S$ of $A$, denoted $\operatorname{Sink}\left(A\left|\, \vec{R},\vec{C} \right.\right)$, is given by $S = D_1 A D_2$, where $D_1 = \operatorname{diag}(r_1, \dots, r_n)$ and $D_2 = \operatorname{diag}(c_1, \dots, c_m)$ are positive diagonal matrices such that $S \vec{1}_m = \vec{R}$ and $S^\top \vec{1}_n = \vec{C}$. When $\vec{R}=\vec{C}=\vec{1}$, we recover the regular Sinkhorn limit $\operatorname{Sink}\left(A\left|\, \vec{1},\vec{1} \right.\right)=\operatorname{Sink}(A)$.
        
    \end{definition}

    Let $S=\operatorname{Sink}\left(A\left|\, \vec{R},\vec{C} \right.\right)$ be given by
    $$
    S=
    \begin{bmatrix}
        s_{11} & s_{12} & \cdots & s_{1m} \\
        s_{21} & s_{22} & \cdots & s_{2m} \\
        \vdots & \vdots & \ddots & \vdots \\
        s_{n1} & s_{n2} & \cdots & s_{nm}
    \end{bmatrix}=
    \begin{bmatrix}
        r_1 a_{11} c_1 & r_1 a_{12} c_2 & \cdots & r_1 a_{1m} c_m \\
        r_2 a_{21} c_1 & r_2 a_{22} c_2 & \cdots & r_2 a_{2m} c_m \\
        \vdots & \vdots & \ddots & \vdots \\
        r_n a_{n1} c_1 & r_n a_{n2} c_2 & \cdots & r_n a_{nm} c_m \\
    \end{bmatrix}
    =D_1 AD_2
    $$

    This matrix equation $S=D_1 AD_2$ for $i = 1, 2, \dots, n \text{ and } j = 1, 2,\dots, m$ gives rise to a set of $n\cdot m$ nonlinear equations,
    \[s_{ij} = r_i a_{ij} c_j \tag{2.1.1}\]\label{2.1.1}

    Additionally, by definition, $\vec{R} = S \vec{1}_m$ and $\vec{C} = S^\top \vec{1}_n$ yield an additional $n+m$ linear equations, where for $i=1,2,\dots,n$ and $j=1,2,\dots, m$,
    \[\sum_{j=1}^m s_{ij} = R_i,\quad \sum_{i=1}^n s_{ij} = C_j \tag{2.1.2}\]\label{2.1.2}
    
    Corresponding equations of \hyperref[2.1.1]{(2.1.1)} can be substituted into \hyperref[2.1.2]{(2.1.2)}, and again for $i=1,2,\dots,n$ and $j=1,2,\dots, m$, this leads to $n+m$ nonlinear equations of the form 
    \[\sum_{j=1}^m r_i a_{ij} c_j = R_i,\quad \sum_{i=1}^n r_i a_{ij} c_j = C_j \tag{2.1.3}\]\label{2.1.3}
    where we are to solve for $r_i$ and $c_j$ in terms of $a_{ij}$, $R_i$, and $C_j$, yielding the closed form of $s_{ij}$ we seek via \hyperref[2.1.1]{(2.1.1)}.\\
    
    Here, we note two things. First, a solution only exists when the consistency condition $\operatorname{sum}\left(\vec{R}\right) = \operatorname{sum}\left(\vec{C}\right)$
    is met, or else the resulting matrix $S$ does not exist. Second, as mentioned prior, Sinkhorn showed that the fact $\operatorname{Sink}\left(A\left|\, \vec{R},\vec{C} \right.\right)=D_1 AD_2$ for positive diagonal matrices $D_{1,2}$ is unique \textit{only up to a scalar multiple}. This is because \hyperref[2.1.3]{system (2.1.3)} is homogeneous in degree 1 for both $r_i$ and $c_j$, or to wit, invariant under the transformation
    $$ r_i\mapsto \lambda r_i,\quad c_j\mapsto \lambda^{-1}c_j,\quad \forall \lambda\in\mathbb R_{>0} $$
    Indeed, we can see how 
    $$\operatorname{Sink}\left(A\left|\, \vec{R},\vec{C} \right.\right) = D_1AD_2 = \left(\lambda D_1\right)A\left(\lambda^{-1}D_2\right) \iff s_{ij}=r_ia_{ij}c_j=\left(\lambda r_i\right)a_{ij}\left(\lambda^{-1}c_j\right)$$

    This is effectively a gauge symmetry of \hyperref[2.1.3]{system (2.1.3)} under the multiplicative group $\mathbb R_{>0}$, and a gauge freedom of the solution set $r_i$, $c_j$. We can thus either 
    \begin{enumerate}
        \item solve for all other variables in terms of some arbitrarily chosen $r_i$ or $c_j$, or
        \item apply a gauge fix, such as designating some arbitrarily chosen $r_i$ or $c_j$ to be a positive real number, ideally as $1$, and solving for all other variables with this gauge fix.
    \end{enumerate}

    Method 2 is typically the easier and more optimal method, as it reduces a lot of algebra.\\

    One last thing to note is that the generalized Sinkhorn limit is symmetric under transposition, and the limits for $n\times m$ matrices are mathematically equivalent to $m\times n$ matrices, shown as follows.
    \begin{theorem}[Proposition 2.1.1 \normalfont{(Generalized Sinkhorn Limit Transpose Invariance)}] \label{Proposition 2.1.1}
        
        Let $A \in \mathbb{R}^{n \times m}_{>0}$ be a positive $n\times m$ rectangular matrix, and let target row and column sums be $\vec{R} \in \mathbb{R}^n_{>0}$ and $\vec{C} \in \mathbb{R}^m_{>0}$, respectively. Let $S = \operatorname{Sink}\left(A\left|\, \vec{R},\vec{C} \right.\right) = D_1 A D_2$ be the generalized Sinkhorn limit of $A$, where $D_1 = \operatorname{diag}(r_1, \dots, r_n)$ and $D_2 = \operatorname{diag}(c_1, \dots, c_m)$ are positive diagonal matrices such that $S \vec{1}_m = \vec{R}$ and $S^\top \vec{1}_n = \vec{C}$.\\
        
        Then, the transposed matrix $A^\top \in \mathbb{R}^{m \times n}_{>0}$ with target row sums $\vec{C}$ and column sums $\vec{R}$ has generalized Sinkhorn limit $\operatorname{Sink}\left(A^\top\left|\, \vec{C},\vec{R} \right.\right) = \operatorname{Sink}\left(A\left|\, \vec{R},\vec{C} \right.\right)^\top = S^\top$.\\

        \textit{Proof}:\quad  We note $A^\top$ is a $m\times n$ matrix, where $n\mapsto m$, $m\mapsto n$. Expanding both sides of the proposed equality gives us by definition
        \vspace*{-1em}
        \begin{align*}
            \operatorname{Sink}\left(A\left|\, \vec{R},\vec{C} \right.\right)^\top &= \operatorname{Sink}\left(A^\top\left|\, \vec{C},\vec{R} \right.\right)\\
            \left(D_1 A D_2\right)^\top &= D_xA^\top D_y\\
            D_2A^\top D_1 &= D_xA^\top D_y
        \end{align*}
        for some unknown $D_{x,y}$. To prove the equality, it suffices to show that for the generalized Sinkhorn limit of $A^\top$, if target row sums are $\vec{C}$ and column sums $\vec{R}$ (swapped compared to $\operatorname{Sink}\left(A\left|\, \vec{R},\vec{C} \right.\right)$, basically), then $D_{x,y} = D_{2,1}$.\\

        From definition, for $\operatorname{Sink}\left(A\left|\, \vec{R},\vec{C} \right.\right)$, row sums are given by $(D_1AD_2)\vec{1}_m=\vec{R}$ and column sums are given by $\left(D_2A^\top D_1\right)\vec{1}_n=\vec{C}$. For the generalized Sinkhorn limit of $A^\top$, its row sums are given by $\left(D_x A^\top D_y\right)\vec{1}_n$, and column sums are given by $\left(D_x A^\top D_y\right)^\top \vec{1}_m = (D_y A D_x) \vec{1}_m$. When row sums are $\left(D_x A^\top D_y\right)\vec{1}_n = \vec{C}$ and column sums are $(D_y A D_x) \vec{1}_m = \vec{R}$, we observe $D_{x,y} = D_{2,1}$.\\
        QED $\blacksquare$

    \end{theorem}

\subsection{Trivial Cases}

    In general, when $A$ is $1\times n$ $\forall n\in\mathbb Z^+$, then the generalized Sinkhorn limit of $A$ is trivial. We can see as follows.\\

    From \hyperref[2.1.3]{system (2.1.3)}, a $1\times n$ matrix gives
    $$
    \begin{cases}
        r_1a_{11}c_1+r_1a_{12}c_2+\dots+r_1a_{1n}c_n=R_1\\
        r_1a_{11}c_1=C_1\\
        r_1a_{12}c_2=C_2\\
        \quad\vdots\\
        r_1a_{1n}c_n=C_n
    \end{cases}
    $$

    Fix $r_1=1$. $R_1=\displaystyle\sum_{n\ge1}C_n$ from consistency. We have $c_j=\dfrac{C_j}{a_{1j}}$ $\forall j\in\mathbb Z^+$ immediately. Therefore,
    \begin{theorem}[Theorem 2.2.1 \normalfont{(Generalized Sinkhorn Limit of $1\times n$ Matrix)}] \label{Theorem 2.2.1}
    
        Let $A = \begin{bmatrix} a_{11} & a_{12} & \dots & a_{1n} \end{bmatrix}$ be a positive $1\times n$ matrix. Then the generalized Sinkhorn limit of $A$ with target row and column sums of $\vec{R}=\begin{bmatrix} R_1 \end{bmatrix} \in \mathbb{R}_{>0}$ and $\vec{C}=\begin{bmatrix} C_1 \\ C_2\\ \vdots\\ C_n \end{bmatrix} \in \mathbb{R}^n_{>0}$, respectively, is given by
        $$\operatorname{Sink}\left(
        \begin{bmatrix}
            a_{11} & a_{12} & \dots & a_{1n} 
        \end{bmatrix}
        \left|\, 
        \begin{bmatrix}
            R_1
        \end{bmatrix},
        \begin{bmatrix}
            C_1 \\C_2\\ \vdots\\ C_n
        \end{bmatrix}
        \right.\right) = 
        \begin{bmatrix}
            C_1 & C_2& \dots& C_n
        \end{bmatrix}
        $$

        \textit{Proof}:\quad  Just jiggle it a bit and the proof should pop out.\\
        QED $\blacksquare$
        
    \end{theorem}
    
    We note that for $n\times 1$ matrices, the generalized Sinkhorn limit follows trivially from \hyperref[Proposition 2.1.1]{Proposition 2.1.1}.

\subsection{$2\times 2$ Matrix} \label{section 2.3}

    Let $A = \begin{bmatrix} a_{11} & a_{12} \\ a_{21} & a_{22} \end{bmatrix}$ be a positive $2\times 2$ matrix. \hyperref[2.1.3]{System (2.1.3)} gives
    \vspace*{-1em}
    $$
    \begin{cases}
        r_1 a_{11} c_1 + r_1 a_{12} c_2 = R_1\\
        r_2 a_{21} c_1 + r_2 a_{22} c_2 = R_2\\
        r_1 a_{11} c_1 + r_2 a_{21} c_1 = C_1\\
        r_1 a_{12} c_2 + r_2 a_{22} c_2 = C_2
    \end{cases}
    $$

    Fix $c_2 = 1$. We thus have 
    $$ 
    \left\{\hspace*{-0.5em}
    \begin{array}{c}
        f_1\\ f_2\\ f_3\\ f_4
    \end{array}\hspace*{-0.5em}
    \right\}=
    \begin{cases}
        r_1 a_{11} c_1 + r_1 a_{12} - R_1\\
        r_2 a_{21} c_1 + r_2 a_{22} - R_2\\
        r_1 a_{11} c_1 + r_2 a_{21} c_1 - C_1\\
        r_1 a_{12} + r_2 a_{22} - C_2
    \end{cases}=0
    $$

    This is a set $\mathcal P = \{f_1,f_2,f_3,f_4\}$ of $4$ highly nontrivial quadratic equations in the polynomial ring $\mathbb R[r_1,r_2,c_1]$. Solving them directly is highly difficult, and even Mathematica's \texttt{Reduce} doesn't output a very nice solution.\\
    
    However, recall that $\mathcal P$ is a generating set for the ideal $I=\langle f_1,f_2,f_3,f_4\rangle \subset\mathbb R[r_1,r_2,c_1]$. For some lexicographic ordering (like $r_1,r_2,c_1$), computing the Gröbner basis $G$ of $I$ gives a much nicer set of polynomials defining the same zero locus as $\mathcal P$, and such $G$ is easily computed in Mathematica (see \hyperref[Mathematica A]{Appendix 4.1.1}).\\

    Furthermore, the choice of lexicographic ordering for the $2\times 2$ case is thankfully not important, as the polynomials are highly symmetric (as is in general for square matrices). This is not the case for a general rectangular matrix, and manual testing may be needed to determine not only the most optimal ordering, but also the most optimal gauge fix.\\

    The computed $G$ contains $8$ polynomials,
    $$
    G=\begin{cases}
        -C_1 - C_2 + R_1 + R_2\\
        -a_{12} a_{22} C_1 - a_{12} a_{21} c_1 C_1 + a_{11} a_{22} c_1 C_2 + a_{11} a_{21} c_1^2 C_2 + a_{12} a_{21} c_1 R_2 - a_{11} a_{22} c_1 R_2\\
        a_{12} a_{21} C_1 - a_{11} a_{21} c_1 C_2 + a_{12} a_{21} a_{22} r_2 - a_{11} a_{22}^2 r_2 - a_{12} a_{21} R_2 + a_{11} a_{22} R_2\\
        a_{22} r_2 + a_{21} c_1 r_2 - R_2\\
        a_{12} C_1 - a_{11} c_1 C_2 + a_{12} a_{22} r_2 + a_{11} a_{22} c_1 r_2 - a_{12} R_2\\
        -C_2 + a_{12} r_1 + a_{22} r_2\\
        -C_1 + a_{11} c_1 r_1 - a_{22} r_2 + R_2\\
        a_{21} C_1 r_2 + a_{11} a_{22} r_1 r_2 + a_{21} a_{22} r_2^2 - a_{11} r_1 R_2 - a_{21} r_2 R_2
    \end{cases}=0
    $$

    Notice that the first polynomial $-C_1 - C_2 + R_1 + R_2=0$ is the consistency condition.\\

    Since we computed $G$ in some lexicographic ordering, solving by starting from the very last polynomial $a_{21} C_1 r_2 + a_{11} a_{22} r_1 r_2 + a_{21} a_{22} r_2^2 - a_{11} r_1 R_2 - a_{21} r_2 R_2=0$ is typically easiest. We solve for $r_1$ (see \hyperref[Mathematica B]{Appendix 4.1.2}), giving us
    \[r_1 = \frac{a_{21}r_2}{a_{11}}\left(\frac{C_1}{R_2-a_{22}r_2}-1\right) \tag{2.3.1} \]

    Then moving onto the second-to-last ($7$th) polynomial, we can substitute in our value of $r_1$ and solve for $c_1$ (see \hyperref[Mathematica C]{Appendix 4.1.3}), giving 
    \[c_1= \frac{R_2-a_{22}r_2}{a_{21} r_2} \tag{2.3.2} \]

    Moving onto the $6$th polynomial, we substitute and solve for $r_2$ (see \hyperref[Mathematica D]{Appendix 4.1.4}), yielding two solutions
    $$
    r_2= \frac{a_{11} a_{22} (C_2+R_2)+a_{12} a_{21} (C_1-R_2)\mp\sqrt{(a_{11} a_{22} (C_2+R_2)+a_{12} a_{21} (C_1-R_2))^2+4 a_{11} a_{22} C_2 R_2 (a_{12} a_{21}-a_{11} a_{22})}}{2 a_{22} (a_{11} a_{22}-a_{12} a_{21})}
    $$

    One of these results is extraneous as it will be negative for some particular choice of positive matrix $A$, so we can plug test matrices until we're able to isolate the extraneous solution. Using the matrix $A=\begin{bmatrix}
        1&2\\
        3&4
    \end{bmatrix}$ and scaling to the regular Sinkhorn limit, we find that for the first (minus square root) result, $r_2\approx 0.11$, and for the second (plus square root) result, $r_2\approx -1.11$, and thus must be extraneous (see \hyperref[Mathematica E]{Appendix 4.1.5}). \\

    Hence we have \[r_2= \frac{a_{11} a_{22} (C_2+R_2)+a_{12} a_{21} (C_1-R_2)-\sqrt{(a_{11} a_{22} (C_2+R_2)+a_{12} a_{21} (C_1-R_2))^2+4 a_{11} a_{22} C_2 R_2 (a_{12} a_{21}-a_{11} a_{22})}}{2 a_{22} (a_{11} a_{22}-a_{12} a_{21})}\tag{2.3.3}\]

    With these results, we solve for $s_{ij}$ (see \hyperref[Mathematica F]{Appendix 4.1.6}) and obtain
    \begin{theorem}[Theorem 2.3.1 \normalfont{(Generalized Sinkhorn Limit of Nonsingular $2\times 2$ Matrix)}] \label{Theorem 2.3.1}
    
        Let $A = \begin{bmatrix} a_{11} & a_{12}\\ a_{21} & a_{22} \end{bmatrix}$ be a positive \textit{nonsingular} $2\times 2$ matrix. Then the generalized Sinkhorn limit $S$ of $A$ with target row and column sums of $\vec{R}=\begin{bmatrix} R_1\\ R_2 \end{bmatrix} \in \mathbb{R}^2_{>0}$ and $\vec{C}=\begin{bmatrix} C_1 \\ C_2 \end{bmatrix} \in \mathbb{R}^2_{>0}$, respectively, is given by
        \begin{align*}S&=
        \operatorname{Sink}\left(
        \begin{bmatrix}
            a_{11} & a_{12}\\
            a_{21} & a_{22}
        \end{bmatrix}
        \left|\, 
        \begin{bmatrix}
            R_1 \\ R_2
        \end{bmatrix},
        \begin{bmatrix}
            C_1 \\ C_2
        \end{bmatrix}
        \right.\right) \\
        &= 
        \begin{bmatrix}
            \dfrac{\alpha(R_2 - 2 C_1 - C_2)+\beta(- R_2 + C_1)+\sqrt{\Delta}}{-2\det(A)} & \dfrac{\alpha(R_2 - C_2)+\beta(-R_2 + C_1 + 2 C_2)-\sqrt{\Delta}}{-2\det(A)} \\
            \dfrac{\alpha(-R_2+C_2)+\beta(R_2+C_1)-\sqrt{\Delta}}{-2\det(A)} & \dfrac{\alpha(-R_2-C_2)+\beta(R_2-C_1)+\sqrt{\Delta}}{-2\det(A)}
        \end{bmatrix}
        \end{align*}
        where $\Delta = (\alpha (C_2+R_2)+\beta (C_1-R_2))^2- 4 \alpha C_2 R_2 \det(A)$, $\alpha = a_{11}a_{22}$, and $\beta = a_{12}a_{21}$.\\

        \textit{Proof}:\quad  Prior in this section.\\
        QED $\blacksquare$
        
    \end{theorem}

    This formula in \hyperref[Theorem 2.3.1]{Theorem 2.3.1} is fully written out in \hyperref[2x2FullMatrixFormula]{Appendix 4.1.7}.\\
    
    As expected, when $\vec{C}=\vec{R}=\vec{1}_2$, we recover Nathanson's $2\times 2$ Sinkhorn limit formula (see \hyperref[Mathematica G]{Appendix 4.1.8}).\\

    However, notice that in \hyperref[Theorem 2.3.1]{Theorem 2.3.1}, $A$ must be nonsingular due to the $\det(A)$ in the denominator. We can further show that when $A$ is singular, the limit of the formula to $A$ gives the generalized Sinkhorn limit of $A$.
    \begin{theorem}[Corollary 2.3.1 \normalfont{(Generalized Sinkhorn Limit of Singular $2\times 2$ Matrix)}]

        Let $A_0 = \begin{bmatrix} a_{11} & a_{12}\\ a_{21} & a_{22} \end{bmatrix}$ be a positive singular matrix (where $\det(A_0)=0$). WLOG choose $a_{12}$ to be a continuous function such that $A(\varepsilon) = \begin{bmatrix} a_{11} & a_{12}(\varepsilon)\\ a_{21} & a_{22} \end{bmatrix}\in \mathbb R^{2\times 2}_{>0}$ for $\varepsilon>0$ is a family of positive nonsingular matrices where $\displaystyle\lim_{\varepsilon\to 0}a_{12}(\varepsilon) = a_{12}$, $\displaystyle\lim_{\varepsilon\to 0} A(\varepsilon) = A_0$, and $\displaystyle\lim_{\varepsilon\to 0} \det(A(\varepsilon)) = 0$. \\
        
        Let target row and column sums be $\vec{R}=\begin{bmatrix} R_1\\ R_2 \end{bmatrix}\in \mathbb{R}^2_{>0}$ and $\vec{C}=\begin{bmatrix} C_1 \\ C_2 \end{bmatrix}\in \mathbb{R}^2_{>0}$, respectively. Then, the analytical continuation of the formula in \hyperref[Theorem 2.3.1]{Theorem 2.3.1} for singular $A_0$ is given by 
        \begin{align*}
            S &= \operatorname{Sink}\left(A_0 \left| \vec{R},\vec{C}\right.\right)\\
            &= \lim_{\varepsilon\to 0}
            \begin{bmatrix}
            \dfrac{\alpha(R_2 - 2 C_1 - C_2)+\beta(\varepsilon)(- R_2 + C_1)+\sqrt{\Delta}}{-2\det(A(\varepsilon))} & \dfrac{\alpha(R_2 - C_2)+\beta(\varepsilon)(-R_2 + C_1 + 2 C_2)-\sqrt{\Delta}}{-2\det(A(\varepsilon))} \\
            \dfrac{\alpha(-R_2+C_2)+\beta(\varepsilon)(R_2+C_1)-\sqrt{\Delta}}{-2\det(A(\varepsilon))} & \dfrac{\alpha(-R_2-C_2)+\beta(\varepsilon)(R_2-C_1)+\sqrt{\Delta}}{-2\det(A(\varepsilon))}
        \end{bmatrix}\\
        &=     \begin{bmatrix}
        C_1-\dfrac{C_1 R_2}{C_1+C_2} & C_2-\dfrac{C_2 R_2}{C_1+C_2} \\
        \dfrac{C_1 R_2}{C_1+C_2} & \dfrac{C_2 R_2}{C_1+C_2} \\
    \end{bmatrix}
        \end{align*}
        where $\Delta = (\alpha (C_2+R_2)+\beta(\varepsilon)(C_1-R_2))^2- 4 \alpha C_2 R_2 \det(A(\varepsilon))$, $\alpha = a_{11}a_{22}$, and $\beta(\varepsilon) = a_{12}(\varepsilon)a_{21}$. The formula in the third line was computed via L'Hôpital's rule (see \hyperref[Mathematica H]{Appendix 4.1.9}).\\

        \textit{Proof}:\quad  For $\Delta = (\alpha (C_2+R_2)+\beta(\varepsilon)(C_1-R_2))^2- 4 \alpha C_2 R_2 \det(A(\varepsilon))$, $\alpha = a_{11}a_{22}$, and $\beta(\varepsilon) = a_{12}(\varepsilon)a_{21}$, let 
        $$S(\varepsilon)=\begin{bmatrix}
            \dfrac{\alpha(R_2 - 2 C_1 - C_2)+\beta(\varepsilon)(- R_2 + C_1)+\sqrt{\Delta}}{-2\det(A(\varepsilon))} & \dfrac{\alpha(R_2 - C_2)+\beta(\varepsilon)(-R_2 + C_1 + 2 C_2)-\sqrt{\Delta}}{-2\det(A(\varepsilon))} \\
            \dfrac{\alpha(-R_2+C_2)+\beta(\varepsilon)(R_2+C_1)-\sqrt{\Delta}}{-2\det(A(\varepsilon))} & \dfrac{\alpha(-R_2-C_2)+\beta(\varepsilon)(R_2-C_1)+\sqrt{\Delta}}{-2\det(A(\varepsilon))}
        \end{bmatrix}$$

        We know that $\forall\varepsilon>0$, $S(\varepsilon) = \operatorname{Sink}\left(A(\varepsilon) \left| \vec{R},\vec{C}\right.\right)$ by \hyperref[Theorem 2.3.1]{Theorem 2.3.1}. We also know that $\operatorname{Sink}$ is a continuous function due to Franklin and Lorenz,\cite{franklin-verygeneral-1989scaling} meaning $$\operatorname{Sink}\left(\displaystyle\lim_{\varepsilon\to0} A(\varepsilon) \left| \vec{R},\vec{C}\right.\right) = \displaystyle\lim_{\varepsilon\to0}\operatorname{Sink}\left( A(\varepsilon) \left| \vec{R},\vec{C}\right.\right)$$

        Therefore, $\operatorname{Sink}\left(A_0 \left| \vec{R},\vec{C}\right.\right) = \displaystyle\lim_{\varepsilon\to0}S(\varepsilon)$.\\
        QED $\blacksquare$
        
    \end{theorem}

    Notice that for singular matrices, the generalized Sinkhorn limit depends only on target row and column sums, and not the initial values of the starting matrix. This is expected, because when a matrix is singular, only proportion information is preserved between entries, and the unique solution is dictated only by the external constraints.\\
    
    We can hence say that in the context of, say, information theory, when the starting matrix is singular, the generalized Sinkhorn limit gives rise to the maximum entropy solution (a matrix that satisfies the marginal constraints but has no extra “shape” imposed by the initial values).\\

    Also I'm... just not gonna bother with $2\times 3$ its way too messy bruh.

\section{Conclusion}
\subsection{Results}

    Through \hyperref[Theorem 2.2.1]{Theorem 2.2.1} and \hyperref[Theorem 2.3.1]{2.3.1}, we have obtained the closed form of the generalized Sinkhorn limit of $1\times n$ and $2\times 2$ matrices. Using the same strategy of utilizing Gr\"obner bases, the closed form of all variables $c_{1,2,3}$ and $r_{2}$ (in this exact lexicographic order) of $2\times 3$ matrices after fixing $r_1=1$ were also computed, though the general closed form matrix was not. Unsurprisingly, the closed form of $r_{2}$ is, for the lack of better words, beyond disgusting, and very large. \\
    
    Furthermore, the simplification process seems highly nontrivial and I am unsure if my result is the fully simplified form. On a machine with 64 GB of RAM, \texttt{FullSimplify} of the overall expression with assumptions uses \textit{\textbf{all}} of the RAM, and stalls after 6-12 hours of running \textbf{\textit{without producing output}}. My result hence has only been run through \texttt{Simplify} with assumptions, and through \texttt{FullSimplify} with assumptions on parts of the expression. For those curious, it is included in the attached Mathematica \texttt{.nb} file.\\

    A Gr\"obner basis was also computed for $2\times 4$ matrices, but due to arbitrary choice in gauge fixing (fixing $r_n$ typically gives best results), lexicographic ordering of variables (ordering $r_n$ last seems to give best results), and choice of Gr\"obner polynomial ordering to solve in (pure trial-and-error), the most optimal path towards a solution is difficult to elucidate, and manual testing probably cannot find said path. As the computation of the Gr\"obner basis follows directly from \hyperref[section 2.3]{section 2.3}, it is omitted. Attempts to solve the system was done through trial-and-error, and thus incredibly messy, without structure, and unrecorded (with many solutions well exceeding 10 MB in size).\\

    However, overall, we managed to "compute" generalized Sinkhorn limits for several small matrices, \textit{in the sense} that we expressed all variables $r_i$, $c_j$ in terms of each other and the constants $a_{ij}, R_i, C_j$, even if not in fully simplified form. The last variable to be solved for is algebraic over the field generated by $a_{ij}, R_i, C_j$ with degree consistent with the binomial pattern $\binom{n+m-2}{n-1}$, and is summarized below.
    \begin{center}\begin{tabular}{|c|c|c|c|}
        \hline
        Matrix Size & Variables & Degree & $\binom{n+m-2}{n-1}$ \\
        \hline
        $1 \times n$ & $n+1$ & $1$ & $1$ \\
        $2 \times 2$ & 4 & 2 & $2$ \\
        $2 \times 3$ & 5 & 3 & $3$ \\
        $2 \times 4$ & 6 & 4 & $4$ \\
        \hline
    \end{tabular}\end{center}

    The degree of the last solved variable also matches the degree of all other variables, since they are all rational functions of each other. This follows from the fact that when computing a Gröbner basis of a $0$ dimensional ideal, all variables belong to the same algebraic extension over some base field. Thus, all entries of a generalized Sinkhorn limit (and by extension, Sinkhorn limits) being combinations of said variables must be algebraic with degree of the last solved variable. This is conjectured by Rowland and Wu in their "Conjecture 33",\cite{rowland-wu-2024} which we now formalize and prove.

    \begin{theorem}[Theorem 3.1.1 \normalfont{(Algebraic Degree of the Generalized Sinkhorn Limit)}] \label{Theorem 3.1.1}

        Let $A \in \mathbb{R}^{n \times m}_{>0}$ be a positive $n\times m$ rectangular matrix $\forall m,n\in\mathbb Z^+$, and let target row and column sums be $\vec{R} \in \mathbb{R}^n_{>0}$ and $\vec{C} \in \mathbb{R}^m_{>0}$, respectively. Let $S = \operatorname{Sink}\left(A\left|\, \vec{R},\vec{C} \right.\right) = D_1 A D_2$ be the generalized Sinkhorn limit of $A$, where $D_1 = \operatorname{diag}(r_1, \dots, r_n)$ and $D_2 = \operatorname{diag}(c_1, \dots, c_m)$ are positive diagonal matrices such that $S \vec{1}_m = \vec{R}$ and $S^\top \vec{1}_n = \vec{C}$.\\
        
        Then, every element $s_{ij}\in S$ is algebraic over the field $\mathbb Q\left(A, \vec{R}, \vec{C}\right)$ with degree at most $\binom{n+m-2}{n-1}$.\\

        \textit{Proof}:\quad  We begin by reinterpreting the problem in projective space, showing that the entries of the limit correspond to points on the intersection of a Segre variety cone and a projective linear subspace. By computing the degree of this intersection, we bound the algebraic degree of each entry of the limit.\\
        
        Since $D_{1,2}$ are diagonal matrices, interpreted as row and column scaling, it is equivalent to $$D_1 A D_2 = \left(\vec{v}_1\vec{v}_2^\top\right)\odot A = V\odot A$$ where $\odot$ is the element-wise product, and $V=\vec{v}_1\vec{v}_2^\top$ so that $S = V\odot A$. Each element of $V$ is $v_{ij} = r_i c_j$, and since $\left(\vec{v}_1\vec{v}_2^\top\right) = \vec{v}_1\otimes \vec{v}_2$ is an outer product, $V$ is rank 1 (all $2\times 2$ minors are null) by construction. Rewriting and homogenizing \hyperref[2.1.3]{system (2.1.3)} with parameter $t$ yields
        $$\left(\sum_{j=1}^m a_{ij} v_{ij}\right) - R_it = 0, \quad \left(\sum_{i=1}^n a_{ij}  v_{ij}\right) - C_j t = 0$$
        These are $n+m$ equations of the new $nm+1$ variables $v_{ij}, t$ whose solutions are equivalence classes $[V: t]\in\mathbb P^{nm}$. However due to the consistency condition, only $n+m-1$ equations are independent.\\ 
        
        Now recall that a Segre embedding $\sigma : \mathbb{P}^n \times \mathbb{P}^m \to \mathbb{P}^{(n+1)(m+1)-1}$
        is given by
        \begin{align*}
            \sigma &: ([X_0 : X_1 : \dots : X_n],[Y_0 : Y_1 :\dots : Y_m ]) \mapsto [X_0Y_0 : X_0Y_1: \dots : X_iY_j: \dots : X_nY_m]\\
            \sigma &: ([X],[Y]) \mapsto [XY]
        \end{align*}
        where we take a pair of equivalence classes in smaller projective spaces and embed them into a larger one. Defining new variables $Z_{ij} = X_i Y_j$ and rewriting the output of a Segre embedding into matrix form, we have
        $$[XY]=[Z] = \begin{bmatrix} Z_{00} & Z_{01} & \dots & Z_{0m} \\ Z_{10} & Z_{11} & \dots & Z_{1m} \\ \vdots & \vdots & \ddots & \vdots \\ Z_{n0} & Z_{n1} & \dots & Z_{nm} \end{bmatrix} = \begin{bmatrix} X_0 \\ X_1 \\ \vdots \\ X_n \end{bmatrix} \begin{bmatrix} Y_0 & Y_1 & \dots & Y_m \end{bmatrix}$$
        The image $\Sigma_{n,m}\subset \mathbb P^{(n+1)(m+1)-1}$ of this embedding is called the Segre variety, and consists of all such equivalence classes of matrices $[Z]$. A natural embedding of the Segre variety into $\mathbb P^{((n+1)(m+1)-1)+1} = \mathbb P^{(n+1)(m+1)}$ is the projective cone $C\left(\Sigma_{n,m}\right)$ over the variety, which consists of all $[Z : \lambda]$ for $[Z] \in \Sigma_{n,m}$. Obviously, the cone has dimension $\dim\left(\Sigma_{n,m}\right)+1 = n+m+1$.\\

        Notice how solutions $[V : t]\in\mathbb P^{nm}$ to our homogenized system lie on the projective cone of some Segre variety $C\left(\Sigma_{n-1,m-1}\right)\in\mathbb P^{nm}$. In particular, we have $n+m-1$ independent equations defining a projective linear subspace $L\subset \mathbb P^{nm}$ of codimension $n+m-1$, to which solutions $[V : t]$ are given by the scheme defined by intersection $H=C\left(\Sigma_{n-1,m-1}\right)\cap L$. This scheme $H$ has dimension
        \vspace*{-0.5em}
        \begin{align*}
            \dim\left(H\right) &= \dim\left(C\left(\textstyle\Sigma_{n-1,m-1}\right)\right) + \dim(L) - \dim\left(\mathbb P^{nm}\right) \\
            &= (n+m-1)+((nm) - (n+m-1)) - nm = 0
        \end{align*}
        as expected. The degree of a Segre variety is defined as the number of intersection points with a linear subspace of complementary codimension in a general position. Moreover, the projective cone has the same degree as the variety, meaning for our variety we have\cite{harris1992algebraic}\cite{MO-179422}
        $$\deg\left(C\left(\textstyle\Sigma_{n-1,m-1}\right)\right) = \deg\left(\textstyle\Sigma_{n-1,m-1}\right) = \binom{(n-1)+(m-1)}{n-1} = \binom{n+m-2}{n-1}$$
        When the subspace is not in a general position, accounting multiplicity returns the same result. \\

        Now, as an example, let $J=\Sigma_{n,m}\cap \mathcal L\in\overline K$ be a 0 dimensional scheme defined over base field $K$, with scheme-theoretic points $P_1,\dots,P_r$ all lying in some common finite extension $K(P)$. Recall that the total number of points $r$ in $J$ with multiplicity, $\deg(J)$, limits how many independent algebraic conditions can be imposed on $P_i$, and thus, bounds the degree of field extension $K(P)/K$. In other words,\cite{MSE-3210719}
        $$[K(P) : K] \le \deg(J) = \deg\left(\Sigma_{n,m}\right) = \binom{n+m}{n}$$

        In our case, for solutions $[V : t] \in H = C\left(\Sigma_{n-1,m-1}\right)\cap L$, our base field is $K = \mathbb Q\left(A, \vec{R}, \vec{C}\right)$ by design. Then, because the cone doesn't change key projective properties (e.g. the Hilbert polynomial, degree, \&c.), the degree of the field extension generated by entries $v_{ij}\in V$ is bounded by $$[K(V):K] \le \binom{n+m-2}{n-1}$$
        QED $\blacksquare$
            
    \end{theorem}

\subsection{Discussion}

    The closed form for the trivial $1\times n$ matrices (and transpose, but this shouldn't need to be said) could be possibly implemented in algorithms as a branch condition to skip computation of these trivial cases, but given that numerical iterative methods (see \hyperref[PythonTestFile]{Appendix 4.2}) are quite efficient even when implemented in "slow" languages like Python, it is unlikely that this implementation will offer much benefits.\\
    
    The closed form for $2\times 2$ matrices derived here serves as a useful pedagogical example, but is impractical for computational purposes. Numerical iterative methods scale matrices much more efficiently than evaluating a complex symbolic formula involving square roots, divisions, and multiplications repeated for each matrix entry. From a computational perspective, the existence of a formula does not outweigh the simplicity and speed of numerical algorithms, \textit{especially} for larger matrices when a formula might not even exist. Thus, the value of the closed form lies more in its theoretical clarity than in its numerical applicability. \\

    We note that once a Gr\"obner basis is found, we can always solve polynomials and substitute variables within the Gr\"obner basis until we arrive at a single univariate polynomial. Then, this polynomial is solved which can be back-substituted into every other variable, and thus obtain a formula for the generalized Sinkhorn limit.\\
    
    As of 2025, clean closed forms of univariate quintics exist in terms of the well understood Bring radical $\operatorname{BR}$, Jacobi theta $\vartheta_n$, or $_4F_3$ hypergeometrics, via the Tschirnhaus transformation. For the sextic, a series solution over $_2F_1$ exists, with potential to transform into a Kamp\'e de F\'eriet $^{p+q}F_{r+s}$ representation.\cite{MSE-Sextic-4728274} For the septic and octic, series solutions over Appell $F_1$ and Lauricella $F_D^{(3)}$, respectively, exist.\cite{MSE-Septic-Octic-4728991}\cite{MO-Overall-Polynomails-450112} This suggests that
    \begin{enumerate}
        \item there exists a closed form (probably with egregiously large algebraic terms) of the generalized Sinkhorn limit for all matrices of the form $1\times n$, $2\times 2$, $2\times 3$, $2\times 4$, and $2\times 5$.
        \item there exists a series representation (again in probably egregiously large algebraic terms) of the generalized Sinkhorn limit of all matrices of the form $2\times 6$, $2\times 7$, $2\times 8$, and $3\times 3$. 
    \end{enumerate}
    For larger matrices, we suggest that the more abstract, polynomial-based framework investigated by Rowland and Wu offer a more viable direction for understanding generalized Sinkhorn limits. Of course, computationally, none of these formulations are practical and a numerical approach is best.\\

    Additionally, in the proof of \hyperref[Theorem 3.1.1]{Theorem 3.1.1}, it was remarked that "when the subspace is not in a general position, accounting multiplicity returns the same result". Given that $A$, $\vec R$, and $\vec{C}$ are all arbitrary positive matrices and vectors, it is likely that some values of them denote some $L$ in a non-general position. Counting multiplicity hides this case, so if we don't account for multiplicity, this seems to suggest that the degree of $V$ is less than the bound. I theorize that this occurs when a matrix is singular, but my brain is absolutely too cooked to flesh this idea out. And I have no clue how $\vec{R}$ or $\vec{C}$ come into play.

\subsection{Acknowledgments}

    I would like to thank \href{https://www.youtube.com/EricRowland}{Rowland's YouTube channel} for introducing me to this problem. I would like to thank members of the Discord servers \textit{Algebraists Anonymous™} and \textit{$>$math is /sci!!!!!!} for providing me with enough motivation to even begin investigating the generalized Sinkhorn limit mathematically, and to prove \hyperref[Theorem 3.1.1]{Theorem 3.1.1} (though admittedly, the proof could be handled much more rigorously and without needing to reference anything).

\section{Appendix}

\subsection{Mathematica $2\times 2$}

\subsubsection{Mathematica A}\label{Mathematica A}

    \lstset{language=Mathematica}
    \begin{lstlisting}
In[1]:= (*List of polynomials to solve*)
        polynomials = {r1*a11*c1 + r1*a12 - row1, r2*a21*c1 + r2*a22 - row2, r1*a11*c1 + r2*a21*c1 - col1, r1*a12 + r2*a22 - col2};

        (*Desired variable elimination order*)
        variables = {r1, r2, c1}; 

        (*Compute the Groebner basis*)
        groebnerResult = GroebnerBasis[polynomials, variables];

        (*Display the Groebner basis*)
        groebnerResult // TableForm

Out[1]//TableForm=
        {{-col1 - col2 + row1 + row2},
         {-a12 a22 col1 - a12 a21 c1 col1 + a11 a22 c1 col2 + a11 a21 c1^2 col2 + a12 a21 c1 row2 - a11 a22 c1 row2},
         {a12 a21 col1 - a11 a21 c1 col2 + a12 a21 a22 r2 - a11 a22^2 r2 - a12 a21 row2 + a11 a22 row2},
         {a22 r2 + a21 c1 r2 - row2},
         {a12 col1 - a11 c1 col2 + a12 a22 r2 + a11 a22 c1 r2 - a12 row2},
         {-col2 + a12 r1 + a22 r2},
         {-col1 + a11 c1 r1 - a22 r2 + row2},
         {a21 col1 r2 + a11 a22 r1 r2 + a21 a22 r2^2 - a11 r1 row2 - a21 r2 row2}}
    \end{lstlisting}

\subsubsection{Mathematica B}\label{Mathematica B}

    \begin{lstlisting}
In[2]:= FullSimplify[
          Solve[
            groebnerResult[[8]] == 0, r1, 
            Assumptions -> {a11 > 0, a12 > 0, a21 > 0, a22 > 0, row1 > 0, row2 > 0, col1 > 0, col2 > 0}
          ]
        ]

Out[2]= {{r1 -> (a21 r2 (-1 + col1/(-a22 r2 + row2)))/a11}}
    \end{lstlisting}

\subsubsection{Mathematica C}\label{Mathematica C}

    \begin{lstlisting}
In[3]:= FullSimplify[
          Solve[
            -col1 + a11 c1 ((a21 r2 (-1 + col1/(-a22 r2 + row2)))/a11) - a22 r2 + row2 == 0, c1, 
            Assumptions -> {a11 > 0, a12 > 0, a21 > 0, a22 > 0, row1 > 0, row2 > 0, col1 > 0, col2 > 0}
          ]
        ]

Out[3]= {{c1 -> (-a22 r2 + row2)/(a21 r2)}}
    \end{lstlisting}

\subsubsection{Mathematica D}\label{Mathematica D}

    \begin{lstlisting}
In[4]:= FullSimplify[
          Solve[
            -col2 + a12 ((a21 r2 (-1 + col1/(-a22 r2 + row2)))/a11) + a22 r2 == 0, r2, 
            Assumptions -> {a11 > 0, a12 > 0, a21 > 0, a22 > 0, row1 > 0, row2 > 0, col1 > 0, col2 > 0}
          ]
        ]

Out[4]= {
         {r2 -> (a12 a21 (col1 - row2) + a11 a22 (col2 + row2) - Sqrt[4 a11 a22 (a12 a21 - a11 a22) col2 row2 + (a12 a21 (col1 - row2) + a11 a22 (col2 + row2))^2])/(2 a22 (-a12 a21 + a11 a22))}, 
         {r2 -> (a12 a21 (col1 - row2) + a11 a22 (col2 + row2) + Sqrt[4 a11 a22 (a12 a21 - a11 a22) col2 row2 + (a12 a21 (col1 - row2) + a11 a22 (col2 + row2))^2])/(2 a22 (-a12 a21 + a11 a22))}
        }
    \end{lstlisting}

\subsubsection{Mathematica E}\label{Mathematica E}

    \begin{lstlisting}
In[5]:= (*Test function for the negative r2 result*)
        testFunctionR2NEGATIVE[a11_, a12_, a21_, a22_, row1_, row2_, col1_, col2_] := (a12 a21 (col1 - row2) + a11 a22 (col2 + row2) - Sqrt[4 a11 a22 (a12 a21 - a11 a22) col2 row2 + (a12 a21 (col1 - row2) + a11 a22 (col2 + row2))^2])/(2 a22 (-a12 a21 + a11 a22));

        (*Test function for the positive r2 result*)
        testFunctionR2POSITIVE[a11_, a12_, a21_, a22_, row1_, row2_, col1_, col2_] := (a12 a21 (col1 - row2) + a11 a22 (col2 + row2) + Sqrt[4 a11 a22 (a12 a21 - a11 a22) col2 row2 + (a12 a21 (col1 - row2) + a11 a22 (col2 + row2))^2])/(2 a22 (-a12 a21 + a11 a22));

        (*Testing the matrix {{1,2},{3,4}} with {row1,row2}={1,1} and {col1,col2}={1,1}*)
        N[testFunctionR2NEGATIVE[1, 2, 3, 4, 1, 1, 1, 1]]
        N[testFunctionR2POSITIVE[1, 2, 3, 4, 1, 1, 1, 1]]
        
Out[5]= 0.112372
Out[6]= -1.11237
    \end{lstlisting}

\subsubsection{Mathematica F}\label{Mathematica F}

    \begin{lstlisting}
In[6]:= (*r2 full definition*)
        r2 = (a12 a21 (col1 - row2) + a11 a22 (col2 + row2) - Sqrt[4 a11 a22 (a12 a21 - a11 a22) col2 row2 + (a12 a21 (col1 - row2) + a11 a22 (col2 + row2))^2])/(2 a22 (-a12 a21 + a11 a22));

        (*sij definition*)
        s11 = FullSimplify[a11 ((a21 r2 (-1 + col1/(-a22 r2 + row2)))/a11) ((-a22 r2 + row2)/(a21 r2)), Assumptions -> {a11 > 0, a12 > 0, a21 > 0, a22 > 0, row1 > 0, row2 > 0, col1 > 0, col2 > 0}];
        s12 = FullSimplify[a12 ((a21 r2 (-1 + col1/(-a22 r2 + row2)))/a11), Assumptions -> {a11 > 0, a12 > 0, a21 > 0, a22 > 0, row1 > 0, row2 > 0, col1 > 0, col2 > 0}];
        s21 = FullSimplify[a21 (r2) ((-a22 r2 + row2)/(a21 r2)), Assumptions -> {a11 > 0, a12 > 0, a21 > 0, a22 > 0, row1 > 0, row2 > 0, col1 > 0, col2 > 0}];
        s22 = FullSimplify[a22 (r2), Assumptions -> {a11 > 0, a12 > 0, a21 > 0, a22 > 0, row1 > 0, row2 > 0, col1 > 0, col2 > 0}];

        (*Find matrix formula*)
        generalizedSinkhornMatrix = {{s11, s12}, {s21, s22}};

        (*Display matrix formula*)
        generalizedSinkhornMatrix // TableForm

Out[7]//TableForm=
        {{
          (a12 a21 (col1 - row2) + a11 a22 (-2 col1 - col2 + row2) + Sqrt[4 a11 a22 (a12 a21 - a11 a22) col2 row2 + (a12 a21 (col1 - row2) + a11 a22 (col2 + row2))^2])/(2 a12 a21 - 2 a11 a22), 
          (a12 a21 (col1 + 2 col2 - row2) + a11 a22 (-col2 + row2) - Sqrt[4 a11 a22 (a12 a21 - a11 a22) col2 row2 + (a12 a21 (col1 - row2) + a11 a22 (col2 + row2))^2])/(2 a12 a21 - 2 a11 a22)},
         {
          (a11 a22 (col2 - row2) + a12 a21 (col1 + row2) - Sqrt[4 a11 a22 (a12 a21 - a11 a22) col2 row2 + (a12 a21 (col1 - row2) + a11 a22 (col2 + row2))^2])/(2 a12 a21 - 2 a11 a22), (a12 a21 (-col1 + row2) - a11 a22 (col2 + row2) + Sqrt[4 a11 a22 (a12 a21 - a11 a22) col2 row2 + (a12 a21 (col1 - row2) + a11 a22 (col2 + row2))^2])/(2 a12 a21 - 2 a11 a22)
        }}
    \end{lstlisting}

\subsubsection{$2\times 2$ Matrix Full Formula}\label{2x2FullMatrixFormula}

    $\operatorname{Sink}\left(A
        \left|\, 
        \vec{R},\vec{C}
        \right.\right)=
        \operatorname{Sink}\left(
        \begin{bmatrix}
            a_{11} & a_{12}\\
            a_{21} & a_{22}
        \end{bmatrix}
        \left|\, 
        \begin{bmatrix}
            R_1 \\ R_2
        \end{bmatrix},
        \begin{bmatrix}
            C_1 \\ C_2
        \end{bmatrix}
        \right.\right) = S = \begin{bmatrix}
        s_{11} & s_{12}\\
        s_{21} & s_{22}
    \end{bmatrix}$ where
    \begin{align*}
        s_{11} &= \dfrac{a_{11}a_{22}(R_2 - 2 C_1 - C_2)+ a_{12}a_{21}(- R_2 + C_1)+\sqrt{(a_{11}a_{22} (C_2+R_2)+ a_{12}a_{21} (C_1-R_2))^2- 4 a_{11}a_{22} C_2 R_2 \det(A)}}{-2\det(A)}\\
        s_{12} &= \dfrac{a_{11}a_{22}(R_2 - C_2)+ a_{12}a_{21}(-R_2 + C_1 + 2 C_2)-\sqrt{(a_{11}a_{22} (C_2+R_2)+ a_{12}a_{21} (C_1-R_2))^2- 4 a_{11}a_{22} C_2 R_2 \det(A)}}{-2\det(A)} \\
        s_{21} &= \dfrac{a_{11}a_{22}(-R_2+C_2)+ a_{12}a_{21}(R_2+C_1)-\sqrt{(a_{11}a_{22} (C_2+R_2)+ a_{12}a_{21} (C_1-R_2))^2- 4 a_{11}a_{22} C_2 R_2 \det(A)}}{-2\det(A)}\\
        s_{22} &= \dfrac{a_{11}a_{22}(-R_2-C_2)+ a_{12}a_{21}(R_2-C_1)+\sqrt{(a_{11}a_{22} (C_2+R_2)+ a_{12}a_{21} (C_1-R_2))^2- 4 a_{11}a_{22} C_2 R_2 \det(A)}}{-2\det(A)}
    \end{align*}

\subsubsection{Mathematica G}\label{Mathematica G}

    \begin{lstlisting}
In[7]:= (*Make a test function for the formula*)
        testFunctionSingular[a11_, a12_, a21_, a22_, row1_, row2_, col1_, col2_] := {{(a12 a21 (col1 - row2) + a11 a22 (-2 col1 - col2 + row2) + Sqrt[4 a11 a22 (a12 a21 - a11 a22) col2 row2 + (a12 a21 (col1 - row2) + a11 a22 (col2 + row2))^2])/(2 a12 a21 - 2 a11 a22), (a12 a21 (col1 + 2 col2 - row2) + a11 a22 (-col2 + row2) - Sqrt[4 a11 a22 (a12 a21 - a11 a22) col2 row2 + (a12 a21 (col1 - row2) + a11 a22 (col2 + row2))^2])/(2 a12 a21 - 2 a11 a22)}, {(a11 a22 (col2 - row2) + a12 a21 (col1 + row2) - Sqrt[4 a11 a22 (a12 a21 - a11 a22) col2 row2 + (a12 a21 (col1 - row2) + a11 a22 (col2 + row2))^2])/(2 a12 a21 - 2 a11 a22), (a12 a21 (-col1 + row2) - a11 a22 (col2 + row2) + Sqrt[4 a11 a22 (a12 a21 - a11 a22) col2 row2 + (a12 a21 (col1 - row2) + a11 a22 (col2 + row2))^2])/(2 a12 a21 - 2 a11 a22)}};

        (*Plug in 1 into the formula, and subtract off Nathanson's formula, the result should be 0*)
        FullSimplify[testFunctionSingular[a11, a12, a21, a22, 1, 1, 1, 1] - 1/(Sqrt[a11 a22] + Sqrt[a12 a21]) {{Sqrt[a11 a22], Sqrt[a12 a21]}, {Sqrt[a12 a21], Sqrt[a11 a22]}}, Assumptions -> {a11 > 0, a12 > 0, a21 > 0, a22 > 0}]

Out[8]= {{0, 0}, {0, 0}}
    \end{lstlisting}

\subsubsection{Mathematica H}\label{Mathematica H}

    \begin{lstlisting}
In[8]:= (*b = a12 a21, a = a11 a22*)
        (*Epsilon -> 0 is the same as b -> a, as both mean the matrix becomes singular*)
        (*We will use b as the active variable*)

        (*s11 case*)
        Limit[(b (col1 - row2) + a (-2 col1 - col2 + row2) + Sqrt[4 a (b - a) col2 row2 + (b (col1 - row2) + a (col2 + row2))^2])/(2 b - 2 a), b -> a]

Out[9]= Indeterminate

In[9]:= (*Thus we can use Lhospital*)
        (*The denominator for all sij will be*)
        FullSimplify[
         D[2 b - 2 a, b],
         Assumptions -> {a > 0, b > 0}
        ]

Out[10]= 2

In[10]:= (*Numerator for s11*)
         FullSimplify[
          D[b (col1 - row2) + a (-2 col1 - col2 + row2) + Sqrt[4 a (b - a) col2 row2 + (b (col1 - row2) + a (col2 + row2))^2], b],
          Assumptions -> {a > 0, b > 0, col1 > 0, col2 > 0, row1 > 0, row2 > 0}
        ]

Out[11]= col1 - row2 + (b (col1 - row2)^2 + a (col2 - row2) row2 + a col1 (col2 + row2))/Sqrt[4 a (-a + b) col2 row2 + (b (col1 - row2) + a (col2 + row2))^2]

In[11]:= FullSimplify[
         Limit[(col1 - row2 + (b (col1 - row2)^2 + a (col2 - row2) row2 + a col1 (col2 + row2))/Sqrt[4 a (-a + b) col2 row2 + (b (col1 - row2) + a (col2 + row2))^2])/2, b -> a],
         Assumptions -> {a > 0, b > 0, col1 > 0, col2 > 0, row1 > 0, row2 > 0}
         ]

Out[12]= col1 - (col1 row2)/(col1 + col2)

In[12]:= (*s12 case*)
         Limit[(b (col1 + 2 col2 - row2) + a (-col2 + row2) - Sqrt[4 a (b - a) col2 row2 + (b (col1 - row2) + a (col2 + row2))^2])/(2 b - 2 a), b -> a]

Out[13]= Indeterminate

In[13]:= (*Thus we can use Lhospital*)
         (*Numerator for s12*)
         FullSimplify[
          D[b (col1 + 2 col2 - row2) + a (-col2 + row2) - Sqrt[4 a (b - a) col2 row2 + (b (col1 - row2) + a (col2 + row2))^2], b],
          Assumptions -> {a > 0, b > 0, col1 > 0, col2 > 0, row1 > 0, row2 > 0}
         ]

Out[14]= col1 + 2 col2 - row2 + (-b (col1 - row2)^2 - a ((col2 - row2) row2 + col1 (col2 + row2)))/Sqrt[4 a (-a + b) col2 row2 + (b (col1 - row2) + a (col2 + row2))^2]

In[14]:= FullSimplify[
         Limit[(col1 + 2 col2 - row2 + (-b (col1 - row2)^2 - a ((col2 - row2) row2 + col1 (col2 + row2)))/Sqrt[4 a (-a + b) col2 row2 + (b (col1 - row2) + a (col2 + row2))^2])/2, b -> a],
         Assumptions -> {a > 0, b > 0, col1 > 0, col2 > 0, row1 > 0, row2 > 0}
         ]

Out[15]= col2 - (col2 row2)/(col1 + col2)

In[15]:= (*s21 case*)
         Limit[(a (col2 - row2) + b (col1 + row2) - Sqrt[4 a (b - a) col2 row2 + (b (col1 - row2) + a (col2 + row2))^2])/(2 b - 2 a), b -> a]

Out[16]= Indeterminate

In[16]:= (*Thus we can use Lhospital*)
         (*Numerator for s21*)
         FullSimplify[
          D[a (col2 - row2) + b (col1 + row2) - Sqrt[4 a (b - a) col2 row2 + (b (col1 - row2) + a (col2 + row2))^2], b],
          Assumptions -> {a > 0, b > 0, col1 > 0, col2 > 0, row1 > 0, row2 > 0}
         ]

Out[17]= col1 + row2 + (-b (col1 - row2)^2 - a ((col2 - row2) row2 + col1 (col2 + row2)))/Sqrt[4 a (-a + b) col2 row2 + (b (col1 - row2) + a (col2 + row2))^2]

In[17]:= FullSimplify[
         Limit[(col1 + row2 + (-b (col1 - row2)^2 - a ((col2 - row2) row2 + col1 (col2 + row2)))/Sqrt[4 a (-a + b) col2 row2 + (b (col1 - row2) + a (col2 + row2))^2])/2, b -> a],
         Assumptions -> {a > 0, b > 0, col1 > 0, col2 > 0, row1 > 0, row2 > 0}
         ]

Out[18]= (col1 row2)/(col1 + col2)

In[18]:= (*s22 case*)
         Limit[(b (-col1 + row2) - a (col2 + row2) + Sqrt[4 a (b - a) col2 row2 + (b (col1 - row2) + a (col2 + row2))^2])/(2 b - 2 a), b -> a]

Out[19]= Indeterminate

In[19]:= (*Thus we can use Lhospital*)
         (*Numerator for s22*)
         FullSimplify[
          D[b (-col1 + row2) - a (col2 + row2) + Sqrt[4 a (b - a) col2 row2 + (b (col1 - row2) + a (col2 + row2))^2], b],
          Assumptions -> {a > 0, b > 0, col1 > 0, col2 > 0, row1 > 0, row2 > 0}
         ]

Out[20]= -col1 + row2 + (b (col1 - row2)^2 + a (col2 - row2) row2 + a col1 (col2 + row2))/Sqrt[4 a (-a + b) col2 row2 + (b (col1 - row2) + a (col2 + row2))^2]

In[20]:= FullSimplify[
         Limit[(-col1 + row2 + (b (col1 - row2)^2 + a (col2 - row2) row2 + a col1 (col2 + row2))/Sqrt[4 a (-a + b) col2 row2 + (b (col1 - row2) + a (col2 + row2))^2])/2, b -> a],
         Assumptions -> {a > 0, b > 0, col1 > 0, col2 > 0, row1 > 0, row2 > 0}
         ]

Out[21]= (col2 row2)/(col1 + col2)

In[21]:= (*Thus we have derived the analytic continuation for singular matrices*)
         singularMatricesSinkhorn2times2[row1_, row2_, col1_, col2_] := {{col1 - (col1 row2)/(col1 + col2), col2 - (col2 row2)/(col1 + col2)}, {(col1 row2)/(col1 + col2), (col2 row2)/(col1 + col2)}};

         (*Here is an example*)
         singularMatricesSinkhorn2times2[30, 30, 10, 50]

Out[22]= {{5, 25}, {5, 25}}
    \end{lstlisting}

\subsection{Python Test File}\label{PythonTestFile}

    This Python script was used to test results numerically. It gives the generalized Sinkhorn limit of any sized matrix.
    \lstset{language=Python}
    \begin{lstlisting}
import numpy as np

def generalized_sinkhorn_knopp(
    matrix: np.ndarray,
    target_row_sums: np.ndarray,
    target_col_sums: np.ndarray,
    tolerance: float = 1e-9,
    max_iterations: int = 1000,
    verbose: bool = True
) -> np.ndarray:
    """
    Calculates the generalized Sinkhorn-Knopp limit for a positive matrix,
    scaling its rows and columns to match specified target sums.

    Iteratively scales rows and then columns until convergence.

    Args:
        matrix (np.ndarray): Initial positive input matrix (2D numpy array).
                             All elements must be strictly positive.
        target_row_sums (np.ndarray): A 1D numpy array of target sums for each row.
                                      Must have the same number of elements as matrix rows.
                                      All elements must be strictly positive.
        target_col_sums (np.ndarray): A 1D numpy array of target sums for each column.
                                      Must have the same number of elements as matrix columns.
                                      All elements must be strictly positive.
        tolerance (float): The convergence criterion. The algorithm stops when the
                           Frobenius norm of the difference between successive matrices
                           is less than this tolerance.
        max_iterations (int): The maximum number of iterations to perform.
        verbose (bool): If True, prints iteration progress.

    Returns:
        np.ndarray: The scaled matrix that approximates the target row and column sums.

    Raises:
        ValueError: If inputs are invalid (e.g., non-positive, incompatible dimensions,
                    or sum mismatch between target row and column sums).
    """

    # --- Input Validation ---
    if not isinstance(matrix, np.ndarray) or matrix.ndim != 2:
        raise ValueError("Input 'matrix' must be a 2D numpy array.")
    if not np.all(matrix > 0):
        raise ValueError(
            "All elements of the input 'matrix' must be strictly positive.")

    if not isinstance(target_row_sums, np.ndarray) or target_row_sums.ndim != 1:
        raise ValueError("'target_row_sums' must be a 1D numpy array.")
    if not np.all(target_row_sums > 0):
        raise ValueError(
            "All elements of 'target_row_sums' must be strictly positive.")

    if not isinstance(target_col_sums, np.ndarray) or target_col_sums.ndim != 1:
        raise ValueError("'target_col_sums' must be a 1D numpy array.")
    if not np.all(target_col_sums > 0):
        raise ValueError(
            "All elements of 'target_col_sums' must be strictly positive.")

    num_rows, num_cols = matrix.shape
    if len(target_row_sums) != num_rows:
        raise ValueError(
            f"Length of 'target_row_sums' ({len(target_row_sums)}) must match "
            f"the number of matrix rows ({num_rows})."
        )
    if len(target_col_sums) != num_cols:
        raise ValueError(
            f"Length of 'target_col_sums' ({len(target_col_sums)}) must match "
            f"the number of matrix columns ({num_cols})."
        )

    # Check the consistency condition (sum of target row sums must equal sum of target column sums)
    if not np.isclose(np.sum(target_row_sums), np.sum(target_col_sums)):
        raise ValueError(
            f"Sum of target row sums ({np.sum(target_row_sums)}) must equal "
            f"sum of target column sums ({np.sum(target_col_sums)}). "
            "A solution does not exist otherwise."
        )

    # --- Algorithm Initialization ---
    current_matrix = np.copy(matrix).astype(
        float)  # Ensure float type for division

    # --- Iterative Scaling ---
    for i in range(max_iterations):
        prev_matrix = np.copy(current_matrix)

        # 1. Row Scaling
        current_row_sums = np.sum(current_matrix, axis=1)
        # Avoid division by zero or very small numbers, though inputs are positive
        # Add a small epsilon to avoid division by zero if sums somehow become zero (just in case)
        row_scaling_factors = target_row_sums / \
            (current_row_sums + np.finfo(float).eps)
        # Apply scaling (multiply each row by its corresponding factor)
        current_matrix = current_matrix * row_scaling_factors[:, np.newaxis]

        # 2. Column Scaling
        current_col_sums = np.sum(current_matrix, axis=0)
        # Add a small epsilon for robustness
        col_scaling_factors = target_col_sums / \
            (current_col_sums + np.finfo(float).eps)
        # Apply scaling (multiply each column by its corresponding factor)
        current_matrix = current_matrix * col_scaling_factors[np.newaxis, :]

        # --- Convergence Check ---
        if np.linalg.norm(current_matrix - prev_matrix, 'fro') < tolerance:
            if verbose:
                print(f"Converged in {i + 1} iterations.")
            break
    else:
        if verbose:
            print(f"Did not converge within {max_iterations} iterations.")

    return current_matrix


# --- Example Usage ---
if __name__ == "__main__":
    # Example 1: Doubly Stochastic (sums to 1)
    print("--- Example 1: Standard Doubly Stochastic (sums to 1) ---")
    A1 = np.array([[2.0, 4], [3, 6.0]])
    R1 = np.array([1.0, 1.0])
    C1 = np.array([1.0, 1.0])

    scaled_A1 = generalized_sinkhorn_knopp(A1, R1, C1)
    print("Original Matrix A1:\n", A1)
    print("Target Row Sums R1:", R1)
    print("Target Column Sums C1:", C1)
    print("\nScaled Matrix A1:\n", scaled_A1)
    print("Resulting Row Sums A1:", np.sum(scaled_A1, axis=1))
    print("Resulting Column Sums A1:", np.sum(scaled_A1, axis=0))
    print("-" * 50)

    # Example 2: Generalized Sinkhorn (arbitrary sums) 
    print("--- Example 2: Generalized Sinkhorn (arbitrary sums) ---")
    A2 = np.array([[1, 2], [400, 9999/17]])
    R2 = np.array([1/3, 16/17])
    C2 = np.array([1/2, 79/102])

    scaled_A2 = generalized_sinkhorn_knopp(A2, R2, C2)
    print("Original Matrix A2:\n", A2)
    print("Target Row Sums R2:", R2)
    print("Target Column Sums C2:", C2)
    print("\nScaled Matrix A2:\n", scaled_A2)
    print("Resulting Row Sums A2:", np.sum(scaled_A2, axis=1))
    print("Resulting Column Sums A2:", np.sum(scaled_A2, axis=0))
    print("-" * 50)
    \end{lstlisting}
    
\printbibliography

\end{document}